\documentclass[11pt]{article}
\usepackage{amsfonts}
\usepackage[mathscr]{eucal}

\newcommand{\n}{\noindent}

\title {$K$-equivalence in Birational Geometry and Characterizations of Complex Elliptic Genera}

\author{
Chin-Lung Wang\thanks{
Department of Mathematics,
National Tsing-Hua University,
Hsinchu, Taiwan.
Email address: dragon@math2.math.nthu.edu.tw
}
}
\date{August, 2000}

\begin{document}
\maketitle

\begin{abstract}
We show that for smooth complex projective varieties the most general combinations
of chern numbers that are invariant under the $K$-equivalence relation consist
of the complex elliptic genera.
\end{abstract}

\bigskip
\centerline{\bf \S1 Introduction and Statements}\medskip

\n The $K$-equivalence relation for $\mathbb{Q}$-Gorenstein projective
varieties was introduced in \cite{Wang} in order to abstract the notion of higher dimensional
{\it composite of flops}. Two birational $\mathbb{Q}$-Gorenstein varieties $X$ and $X'$ are
$K$-equivalent (denoted by $X=_K X'$) if there is a
smooth birational model $\phi:Y\to X$ and $\phi':Y\to X'$ such that
$\phi^*K_X=_\mathbb{Q}\phi'^*K_{X'}$. This simple notion appears naturally in birational
geometry and we are interested in characterizing those geometric/topological
invariants that are invariant under it. \par

The basic strategy formulated in \cite{Wang} is
a {\it meta theorem}: in order to obtain numerical invariants
under the $K$-equivalence relation, it suffices to have a suitable integration theory
which admits a {\it nice} change of variable formula for birational morphisms and
geometric/topological interpretations of the integration. For $K$-equivalent smooth complex
projective varieties, the invariance of Betti and Hodge numbers has been verified in
\cite{Wang} \cite{DeLo} \cite{Batyrev}.
Based on this, we further conjectured that their ordinary cohomology groups
are canonically isomorphic under the cohomology correspondence induced from the graph
of the given birational map. It is clear that this map respects the $\mathbb{Q}$-Hodge structures.\par

Not all interesting invariants are invariant under the
$K$-equivalence relation. For example, while the cohomology groups for smooth
threefolds are canonically isomorphic under a classical flop, their ring structures are in general
different. We find two recent works that may lead to an explanation. One is the
{\it quantum minimal model conjecture} raised by Ruan \cite{Ruan} (cf.\ \S6).
Another one is Totaro's work \cite{Totaro}: {\it complex elliptic genera $=$
complex cobordism ring modulo classical flops}, or equivalently chern numbers
invariant under classical flops consist of precisely the complex elliptic genera.\par

We show in this paper a generalization of Totaro's result: {\it complex elliptic genera $=$ complex cobordism
ring modulo $K$-equivalence}. A surprising conclusion is that in the complex cobordism ring, the ideal
generated by classical flops equals the seemingly much larger ideal generated by
all $K$-equivalence pairs. To summarize our approach, we follow the meta theorem.
The target invariants
are genera (or chern numbers), which by definition are certain integrals, so we only need a
nice change of variable formula. This is treated in two steps: \par

Inspired by the work of Hirzebruch et.\ al.\ \cite{HBJ} on the characterization of the (real)
elliptic genera (of Landweber and Stong) as the most
general genera that are multiplicative on fiber bundles with fiber $\mathbb{P}^{2n -1}$ for
all $n\in\mathbb{N}$, we characterize the complex elliptic genera
(studied by Witten, Hirzebruch and subsequently by Krichever, H\"ohn and Totaro)
in a similar flavor via blowing-ups along smooth centers.\par

Let $X$ be a compact complex manifold or
a proper smooth variety (over field of arbitrary characteristic). For a commutative ring $R$,
an $R$-genus $\varphi$ is defined by a power series $Q(x)\in R[\![x]\!]$ through Hirzebruch's
multiplicative sequence $K_Q$ (or $K_\varphi$). As usual we write
$Q(x) = x/f(x)$. \par\medskip

\n{\bf Theorem A} (Residue Theorem) {\it For any cycle $D$ in $X$ and for any blowing-up
$\phi: Y\to X$ along smooth center $Z$ with exceptional divisor $E$, one has for any power
series $A(t)\in R[\![t]\!]$:
\begin{eqnarray*}
\int_{\phi^*D}A(E)\,K_Q(c(T_Y))
 \!\!&=&\!\! \int_D A(0)\,K_Q(c(T_X))\\
 && + \int_{Z.D} {\rm Res}_{\,t=0}\Big({A(t)\over f(t)\prod_{i=1}^r f(n_i - t)}\Big)
 \,K_Q(c(T_Z)).
\end{eqnarray*}
Here $n_i$'s denote the formal chern roots of the normal
bundle $N_{Z/X}$ and the residue stands for the coefficient of the degree $-1$ term of a
Laurent power series with coefficients in the cohomology ring or the Chow ring of $X$.}\par\medskip

\n{\bf Theorem B} (Characterization of Complex Elliptic Genera) {\it Consider the following
sets of power series $f(x)=x+\cdots\in\mathbb{C}[\![x]\!]$ (or $\mathbb{C}$-genera $\phi$'s):
\begin{itemize}
\item[$S_1:$] $\varphi$ admits a first step change of variable formula.
That is, for each $r\in\mathbb{N}$ there exists a power series
$A(t,r)$ in $t$ serves as the Jacobian factor such that $A(0,r) = 1$ and
$$
\int_X K_\varphi(c(T_X)) = \int_Y A(E, r)\,K_\varphi(c(T_Y))
$$
for any blowing-up $\phi: Y\to X$ along smooth center of codimension $r$ with exceptional divisor $E$.
\item[$S_2:$] $f(x)$ is a solution to the following functional equation:
$$
{1\over f(x)f(y)} = {A(x)\over f(x)f(y - x)} + {A(y)\over f(y)f(x - y)}
$$
for some power series $A(t)$ with $A(0) =1$.
\item[$S_3:$] $\varphi$ is a specialization of the complex elliptic genera. That is,
$\varphi$ is parameterized by $(k,a,b,g_2)\in\mathbb{C}^4$ such that
$$
{f'(x)\over f(x)}=-{1\over 2}
{\wp'(x)-b\over \wp(x)-a} + k.
$$
Where $\wp(x)$ is the unique function with a pole of order 2 at zero such that
$\wp'(x)^2=4\wp(x) -g_2 \wp(x)-g_3$ with $g_3$ defined by $b^2=4a^3 -g_2a-g_3$. Equivalently,
$\varphi$ is parametrized by $k\in\mathbb{C}$ and affine Weierstrass equations, which may define
singular curves, with a marked point.
\end{itemize}
Then $S_1\subset S_2\subset S_3\cong\mathbb{C}^4$. Moreover, $S_2$ contains precisely those
$f$'s with $(a,b)$ not a 2 torsion point plus exceptional cases $e^{kx}\sinh(sx)/s$ with
$s^2=6f_3=f'''(0)$. In the former cases, $A(x)$ is uniquely determined by $f$ and for the
exceptional cases $A(x) = e^{-kx}(a_1\sinh(sx)/s +\cosh(sx))$ where $a_1=A'(0)$. $S_1$ contains all $f$'s
with $(a,b)$ a non-torsion point. In particular, the generic point of the complex elliptic genera
(the power series $f$ with parameters) is a solution to $S_1$ and $S_2$.}\par\medskip

Let $(a,b)=(\wp(z),\wp'(z))$, in terms of the Weierstrass $\sigma$ function we will write down $f$ and
find a candidate for $A(t,r)$ in \S4, namely
$$
f(x)=e^{(k+\zeta(z))x}\,{\sigma(x)\sigma(z)\over\sigma(x+z)},\qquad
A(t,r) = e^{-(r-1)(k+\zeta(z))t}\,{\sigma(t + rz)\sigma(z)\over\sigma(t+z)\sigma(rz)}.
$$
This shows that for blowing-ups along codimension $r$ centers, there could be {\it no}
change of variable formula when one specializes $z$ to an $r$ torsion point.\par

The second step is the {\it full} change of variable formula for birational morphisms. At this moment
we can only prove this for morphisms that can be factorized into
composite of blowing-ups and blowing-downs along smooth centers. In the
algebraic case, we may complete the proof in the characteristic zero case using recent result
of Wlodarsczyk et.\ al.\ \cite{AKMW} \cite{Wlodarsczyk} on the {\it weak factorization theorem}.\par\medskip

\n{\bf Theorem C} (Change of Variable Formula) {\it Let $\varphi$ be the complex elliptic genera. Then
for any algebraic cycle $D$ in $X$ and birational morphism $\phi:Y \to X$ with
$K_Y = \phi^*K_X + \sum e_i E_i$, we have
$$
\int_D K_\varphi(c(T_X))
= \int_{\phi^*D} \prod\nolimits_i A(E_i, e_i + 1)\,K_\varphi(c(T_Y)).
$$
Or equivalently, $\phi_*\prod_i A(E_i, e_i + 1)\,K_\varphi(c(T_Y))=K_\varphi(c(T_X))$.}\par\medskip

It is well-known that the Todd genus is the only complex genus that is absolutely birationally
invariant. In this case, $Q(x) = x/(1 - e^{-x})$, $A(t,r) = 1 = A(x)$ and Theorem C is
a special case of the Grothendieck-Riemann-Roch Formula $\phi_*{\rm td}(T_Y) = {\rm td}(T_X)$,
which is valid without any restriction on the ground field. Our guiding principle is that one
should regard the Todd genus as a {\it rational measure} and the complex elliptic genera as
the {\it elliptic measure} --- think of $K_\varphi(c(T_X))$ as $d\mu_X$!
Theorem C together with Totaro's result \cite{Totaro} imply \par\medskip

\n{\bf Theorem D} (Invariance under $K$-equivalence) {\it The most general
chern numbers that are preserved under the $K$-equivalence relation among smooth proper
complex varieties consist of the complex elliptic genera. Moreover,
complex elliptic genera are precisely those genera that could be defined on log-terminal singular
varieties through the change of variable formula, with the result to be independent of the chosen
smooth model.}\par\medskip

Our original plan is to regard Theorem C as a Grothendieck
Riemann-Roch type problem and prove Theorem D via Theorem B and C only. One reason for doing so
is to get a characteristic free treatment of all these results in the algebraic case. We
realize only part of it in the current work. Also it should be interesting to
compare our approach with Totaro's, which used the rigidity property of the complex
elliptic genera.
\par

Historically, physicists in the 80's has predicted that the
two parameter elliptic genera (with $k=0$) can be
defined on singular Calabi-Yau varieties (at least for orbifolds)
and should agree with the value on the
(non-unique) crepant resolutions, if there are any. This has motivated many earlier works
on this subject. In this context, a recent preprint of Borisov and Libgober
\cite{BoLi} contains results that are related to our present work. In particular, they also
showed that elliptic genera (with $k=0$) could be defined on varieties with at most
log-terminal singularities.\par\medskip

\n{\bf Acknowledgement} I would like to thank Professors Esnault and Viehweg for bringing
Totaro's work into my attention while I was visiting them at the University of Essen
in January 2000.
Part of this work was done there. Also I would like to thank J.-K. Yu for extremely helpful
discussions on the functional equations considered in \S3. Especially the given proof of
$S_2\subset S_3$ in Theorem B is a joint work with him, which is better then my
original approach. The National Center of Theoretic
Sciences (Hsinchu, Taiwan) has provided wonderful working conditions during the summer
of year 2000 while the main part of this work was completed there.\\

\centerline{\bf \S2 A Residue Theorem for Genera}\medskip

Let $X$ be a compact complex manifold or a proper smooth algebraic variety over an
algebraically closed field.
Let $c(T_X) = \prod\nolimits_i (1 + x_i)$ be the formal decomposition into chern roots.
For any commutative ring $R$, an $R$-genus $\varphi(X)$ is defined by a power series
$Q(x) = 1 + \cdots\in R[\![x]\!]$ through Hirzebruch's recipe of multiplicative sequence
$$
\varphi(X) = \int_X K_\varphi(c(T_X)) = \int_X \prod\nolimits_i Q(x_i).
$$
We will also write $K_\varphi$ as $K_Q$ if $Q$ defines $\varphi$.\par

In order to represent $\varphi(X)$ in terms of data on $Y$ under a blowing-up $\phi:Y\to X$
along $Z$, we will need a localization result. To begin with, let us recall and extend some
results from \cite{HBJ}, Ch.3. As usual, let $Q(x) = x/f(x)$ with $f(x) = x + \cdots$. Then $f(u+v) =F(f(u),f(v))$ for an unique
power series $F(y_1,y_2)=\sum_{(r,s)\ne(0,0)} a_{rs}y_1^ry_2^s = y_1+y_2 +\cdots$.
For cycles $u_j\in H^2(M,\mathbb{Z})$ (or $Z_{n-1}(M)$), the virtual genus is defined by
$\varphi(\prod u_j):=\int_M K_Q(c(T_M))\prod f(u_j)$. When $\prod u_j$ is represented by
a smooth subvariety $V$, $\varphi(\prod u_j)\equiv \varphi(V)$. Then one has \cite{HBJ}, p.37:
$$
\varphi(u + v)=\sum\nolimits_{(r,,s)\ne (0,0)} a_{rs}\varphi(u^rv^s).
$$
Moreover, for $D\subset M$ an analytic/algebraic cycle, we may define
$\varphi_D(\prod u_j)$ by $\int_D K_Q(c(T_M))\prod f(u_j)$. Then the above is also
valid for $\varphi_D$.\par

Now let $g := f^{-1}$, the inverse power series (so $g'(y) = \sum_{n\ge 0}\varphi(\mathbb{P}^n)y^n$).
Let $H_{ij}\subset \mathbb{P}^i\times\mathbb{P}^j$ be the degree $(1,1)$ hypersurfaces,
then by \cite{HBJ}, p.39:
$$
F(y_1,y_2)g'(y_1)g'(y_2) = \sum\nolimits_{(i,j)\ne(0,0)}\varphi(H_{ij})y_1^iy_2^j.
$$

\n{\bf Proposition 2.1} {\it For any complex genus $\varphi$ and two blowing-ups
$\phi:Y\to X$ and $\phi':Y'\to X'$ along isomorphic smooth center $Z\cong Z'$,
if $N_{Z/X}\cong N_{Z'/X'}$ and there are cycles $D\subset X$, $D'\subset X'$ such that $Z.D = Z'.D'$,
then $\varphi_D(X) - \varphi_{\phi^*D}(Y) = \varphi_{D'}(X') - \varphi_{\phi'^*D'}(Y')$.}\par\medskip

\n{\it Proof.}
It is clear that one may reduce the proof to the case that $X' = \mathbb{P}_Z(N\oplus{\bf 1})$ and
$Y' = {\rm Bl}_Z \mathbb{P}_Z(N\oplus{\bf 1})$, where $N$ denotes the normal bundle of
$Z$ and ${\bf 1}$ is the trivial bundle. To get relations among these spaces, we apply
deformation to the normal cone to the inclusion $i:Z\to X$ \cite{Fulton}. That is,
we form the blowing-up $\Phi:M\to X\times\mathbb{P}^1$ along $Z\times\{\infty\}$. Then
$X \cong M_0 \sim M_\infty = Y\cup\mathbb{P}_Z(N\oplus{\bf 1}) =: u + v$, with
$uv = \mathbb{P}_Z(N)$ and $(u + v)u = 0 = (u + v)v$. If we plug in $y_1 = f(u)$ and $y_2 = f(v)$
into the previous formula, then since $\varphi((u+v)u)=0=\varphi((u+v)v)$ and
$u^i=(-1)^{i-1}uv^{i-1}$, $v^j=-uv^{j-1}$, we see that only the constant term of $g'$
contributes and then
\begin{eqnarray*}
\varphi(X) \!\!\!&=&\!\! \varphi (u+v)=\int_M f(u+v)\,K_Q(c(T_M))\\
  &=&\!\!\! \int_M f(u+v)g'(f(u))g'(f(v))\,K_Q(c(T_M))\\
  &=&\!\!\! \int_M F(f(u),f(v))g'(f(u))g'(f(v))\,K_Q(c(T_M)) \\
  &=&\!\!\! \sum\nolimits_{(i,j)\ne(0,0)}\varphi(H_{ij})\varphi(u^iv^j)\\
  &=&\!\!\! \varphi(Y) + \varphi(\mathbb{P}_Z(N\oplus{\bf 1}))
        + \sum\nolimits_{i+j\ge 2}(-1)^{i-1}\varphi(H_{ij})\varphi(E^{i+j-1}).
\end{eqnarray*}
Notice that this is almost a change of variable formula except the appearance of the term
$\varphi(\mathbb{P}_Z(N\oplus{\bf 1}))$. If we apply this formula to the inclusion
$i':Z\to N\to X' = \mathbb{P}_Z(N\oplus{\bf 1})$ with $Y' = {\rm Bl}_Z \mathbb{P}_Z(N\oplus{\bf 1})$, then
the last term is exactly $-\varphi(Y')$. Hence in particular
$$
\varphi(X)-\varphi(Y) = \varphi(\mathbb{P}_Z(N\oplus{\bf 1})) -
\varphi({\rm Bl}_Z \mathbb{P}_Z(N\oplus{\bf 1})).
$$

If instead of the total space $M$ we use the cycle $\Phi^*(D\times \mathbb{P}^1)$, by restricting
the above computation to this cycle, we get the result. \hfill QED\\

\n{\bf Remark 2.2} Proposition 2.1 should be a standard fact at least when $D=X$ and $D'=X'$,
and it is quite visualizable in the case of complex manifolds. The above proof
works in both analytic and algebraic cases.\par\medskip

With this proposition, it suffices to prove Theorem A for any other pair $Y'\to X'$.
First of all, if $A(0) = 0$, then the whole thing localizes to $E\cong \mathbb{P}_Z(N)\to Z$
and becomes a fiber integration formula, so it does not depend on which pair $Y\to X$ we use.
If $A(0)\ne 0$, we may normalize it so that $A(0) = 1$ and we may rewrite the power series $A(t)$ as
$J(f(t))$, then the equality
$$
\int_{X'} K_\varphi(c(T_{X'})) = \int_{Y'} J(f(E))\,K_\varphi(c(T_{Y'})) + \int_Z {\rm Res}
$$
and Proposition 2.1 implies that $\varphi(X) - \varphi(Y) = \varphi(X') - \varphi(Y')$ equals
$$
\int_{Y'} (J(f(E)) - 1)\,K_\varphi(c(T_{Y'})) + \int_Z {\rm Res} = \varphi(J(E) - 1)
+ \int_Z {\rm Res}.
$$
All terms in the last expression of virtual genus involves $E^i$ with $i\ge 1$, hence localize to $E$
and equals the same expression $\varphi(J(E)) - 1) + \int_Z {\rm Res}$ but with ambient space $Y$.
Write it back we get the first step change of variable formula we want.
The same argument applies to the general case involving $D$ etc.\ too.\par

So we may assume that $X = \mathbb{P}_Z(N\oplus{\bf 1})$ with $i:Z\to X$
the embedding as the zero section $Z\to N\to \mathbb{P}_Z(N\oplus{\bf 1})$. Then
$N = i^*{\cal Q}$ where ${\cal Q}$ is the universal quotient bundle in the tautological sequence
$$
0\to {\cal S}\to p^*(N\oplus{\bf 1})\to {\cal Q}\to 0.
$$
Here $p:\mathbb{P}_Z(N\oplus{\bf 1})\to Z$ is the projection map. In such case, there
are explicit formula relating $c(T_X)$ and $c(T_Y)$ (\cite{Fulton}; p.301, (a)):\par\medskip

\n{\bf Lemma 2.3} {\it If in a blowing-up $\phi: Y\to X$ along smooth center $Z$, the normal
bundle $N = i^*{\cal Q}$ for some vector bundle ${\cal Q}$ on $X$, then
$$
c(T_Y) = \phi^* c(T_X)\,\phi^*c({\cal Q})^{-1}\,(1 + E)\,c(\phi^*{\cal Q}\otimes\mathscr{O}(-E)),
$$
where $E$ is the exceptional divisor.}\par\medskip

Let $q_i$'s be the formal chern roots of ${\cal Q}$ and let
$$
R(t)=\sum_{k\ge0} R_k(\phi^*{\cal Q})\,t^k = Q(t)\prod\nolimits_{i=1}^r Q(\phi^*q_i - t)A(t),
$$
a power series with coefficients in the cohomology ring $H^*(Y,\mathbb{Q})$ (or in the Chow ring). Then
\begin{eqnarray*}
&&\!\!\!\int_Y  A(E)\,K_Q(c(T_Y))\\
&=&\!\!\! \int_Y \Big[Q(E)\prod\nolimits_{i=1}^r Q(\phi^*q_i - E)A(E)\Big]\,
  \phi^*\Big(K_Q(c(T_X))\,K_Q(c({\cal Q}))^{-1}\Big)\\
&=&\!\!\! \int_Y R(E)\,\phi^*\Big(K_Q(c(T_X))\,K_Q(c({\cal Q}))^{-1}\Big).
\end{eqnarray*}
It is clear that the constant term gives rise to the main term:
$$
\int_Y A(0)\,\phi^*K_Q(c(T_X)) = \int_X A(0)\,K_Q(c(T_X)).
$$
Let $j:E\to Y$ be the inclusion map and $\bar\phi:E=\mathbb{P}_Z(N)\to Z$ be the
restriction of $\phi$. One has $j^*\phi^*{\cal Q} = \bar\phi^*i^*{\cal Q} = \bar\phi^* N$
and also
$$
0\to\bar\phi^*T_Z \to \bar\phi^*i^*T_X\equiv j^*\phi^*T_X\to \bar\phi^*N\to 0,
$$
hence that $j^*\phi^*(K_Q(c(T_X))\,K_Q(c({\cal Q}))^{-1}) = \bar\phi^*K_Q(c(T_Z))$.
If $1$ is the fundamental class of $E$ then
$j_*(1) = E$ and $e = E|_E = j^*E$ is the class of $\mathscr{O}_{\mathbb{P}(N)}(-1)$.
Since $\alpha\cdot E = \alpha\cdot j_*(1) = j^*\alpha$, all
non-constant terms localize to $E$. Let $R_1(t) = \sum_{k\ge 1}R_k(\bar\phi^*N,r)\,t^{k-1}$.
Then the remaining terms give rise to
$$
{\cal R} = \int_E  R_1(e)\,\bar\phi^*K_Q(c(T_Z)) = \int_Z \bar\phi_*R_1(e)\,K_Q(c(T_Z)).
$$

To perform the fiber integration, we make use of Segre classes as in \cite{Totaro}, Lemma 5.1:\par\medskip

\n{\bf Lemma 2.4} (\cite{Fulton}, p.47) {\it Let $s(N) = \sum s_k(N)$ such that $s(N)c(N)=1$, then
\begin{itemize}
\item[\rm (1)] $\bar\phi_*e^k = 0$ for $0\le k\le r-2$ and
\item[\rm (2)] $\bar\phi_*e^{(r-1) + k} = (-1)^{(r-1) + k} s_k(N)$ for $k\ge 0$.
\end{itemize}}

Apply this Lemma we get
$$
{\cal R} = -\int_Z \Big[\sum\nolimits_{k\ge 0}(-1)^{r+k}R_{r+k}\,s_k(N)\Big]\,K_Q(c(T_Z)).
$$
Now $R_{r+k} = R_{r+k}(N)$ is the coefficient of degree $r+k$ term in
$$
R(t) = Q(t)\prod\nolimits_{i=1}^r Q(n_i - t)\,A(t).
$$
And as power series, $s_t(N)c_t(N)=(\sum s_k(N)\,t^k)\prod(1 + n_it)=1$, which
when replacing $t$ by $-1/t$, one gets
$$
\Big(\sum\nolimits_{k\ge 0}(-1)^k s_k(N)\,{1\over t^k}\Big)\prod\Big(1 - {n_i\over t}\Big)=1.
$$
This implies that $\bar\phi_*R_1(e) = -[\cdots]$ is $(-1)^{r + 1}$ times the coefficient of
degree $r$ term of $R(t)S_{-1/t}(N)$, that is, degree $-1$ term of the Laurent power series in $t$:
$$
Q(t)\prod\nolimits_{i=1}^r Q(n_i - t)\,A(t)\,{1\over t^{r+1}}\,{1\over\prod_{i=1}^r(1-{n_i\over t})}
=:{A(t)\over f(t)\prod_{i=1}^r f(n_i - t)}.
$$
Here we agree to denote that
$$
{1\over n_i -t} := {-1\over t}{1\over 1 - n_i/t} := -{1\over t} - {n_i\over t^2} - {n_i^2 \over t^3} - \cdots
\hbox{\rm(finite terms)}.
$$
So the pole still occurs at $t=0$.\par

Notice that in all the above computations, we may replace the integrals on $X$, $Y$, $Z$ and $E$
by $D$, $\phi^*D$, $i^*D = Z.D$ and $E.\phi^*D = \bar\phi^*(i^*D)=\bar\phi^*(Z.D)$.
we have thus proved Theorem A. \hfill QED\\

\centerline{\bf \S3 Functional Equations versus Differential Equations}\medskip

In this and next sections we give the proof of Theorem B. We will show that
$S_1\subset S_2\subset S_3$ in this section.\par\medskip

For $S_1\subset S_2$, let $r\ge 2$ be a fixed integer. By Theorem A, the defining
condition of $S_1$ is equivalent to that, for any proper smooth variety $Z$
(of arbitrary dimension) and any rank $r$ vector bundle $N\to Z$, the coefficient of the degree
$-1$ term (residue at $t=0$) of $\displaystyle{A(t,r)\over f(t)\prod_{i=1}^r f(n_i - t)}$ is zero.
By our definition of $\displaystyle{1\over n_i - t}$, in order to compute the residue at $t=0$, we may
treat $n_i$ as distinct complex parameters and then compute the total residue at $t=0$ and
$t=n_1,\ldots,n_r$. By the Cauchy residue theorem, the total residue is given by
$$
{1\over \prod_{i=1}^r f(n_i)} - \sum\nolimits_{j=1}^r {A(n_j,r)\over f(n_j)\prod_{i\ne j} f(n_i - n_j)}.
$$
So the defining property of $S_1$ is equivalent to that as formal power series in $x_i$'s, we have
the functional equations:
$$
{1\over \prod_{i=1}^r f(x_i)} = \sum\nolimits_{j=1}^r {A(x_j,r)\over f(x_j)\prod_{i\ne j} f(x_i - x_j)}.
$$
When $r=2$, let $A(t) = A(t, 2)$. The vanishing of the residue for arbitrary rank two bundle
$N\to Z$ is then equivalent to the following power series identity:
$$
{1\over f(x)f(y)} = {A(x)\over f(x)f(y - x)} + {A(y)\over f(y)f(x - y)}.
$$
This is just the functional equation defining $S_2$.\par

For $S_2\subset S_3$, let $f(x) = \sum_{i\ge 1} f_i x^i$ with $f_1 = 1$ and
$A(x)=\sum_{i\ge 0} a_i x^i$ with $a_0 = 1$. The functional equation is equivalent to FE:
$$
f(x - y)f(y - x) = A(x)f(y)f(x - y) + A(y)f(x)f(y - x).
$$
To solve all such $f$ and $A$, our guiding principle is that an identity of this type should be
closely related to the addition law of Weierstrass $\wp$ functions. Since elliptic functions
of such type always satisfy certain ordinary differential equations, we will try to transform our
functional equation into certain ODE's so that we can explicitly solve them. Before doing so,
we notice that if $(f(x),A(x))$ is a solution, then $(e^{kx}f(x),e^{-kx}A(x))$ will also be a solution.
Hence without loss of generality we may assume that $f_2 =0$. \par

First of all, by differentiating FE in $y$ and set $y = 0$, one gets FE':
$$
A(x)f(x) + a_1f(-x)f(x) + f(-x)f'(x) = 0,
$$
That is, $A(x)$ is completely determined by $f$ and $a_1$. Now we plug in $A(x)$ and $A(y)$
into FE and get FE*:
\begin{eqnarray*}
 &&f(x)f(y)f(x-y)f(y-x)\\
 &&\hskip 20pt +\big(a_1f(x)f(-x) + f'(x)f(-x)\big)f(x-y)f(y)^2 \\
 &&\hskip 40pt +\big(a_1f(y)f(-y) + f'(y)f(-y)\big)f(y-x)f(x)^2 = 0.
\end{eqnarray*}
Now we differentiate FE* in $y$ twice and set $y=x$ to get DE-1:
\begin{eqnarray*}
1 + a_1\big(f(x)f'(-x) + f'(x)f(-x)\big) \!\!\!\!&+&\!\!\!\! f'(x)f'(-x)\\
 &+&\!\!\!\! \Big(2{f(-x)f'(x)^2\over f(x)} -f(-x)f''(x)\Big) = 0.
\end{eqnarray*}
If we differentiate FE* in $y$ three times and set $y=0$, we get DE-2:
\begin{eqnarray*}
6f_3f(x)f(-x) + 2a_1\big(f(x)f'(-x) \!\!\!\!&+&\!\!\!\! f'(x)f(-x)\big) + 2f'(x)f'(-x) \\
&+&\!\!\! \!\Big(2{f(-x)f'(x)^2\over f(x)} -f(-x)f''(x)\Big) = 0.
\end{eqnarray*}

To motivate the following calculations, notice that $\wp(x)$
is even with principal part $1/x^2$ at $x = 0$ and with no constant terms. Since
$$
{-1\over f(x)f(-x)}={1\over x^2}(1 - f_3x^2  +\cdots)(1 - f_3x^2 + \cdots)={1\over x^2} - 2f_3 + \cdots,
$$
an ambitious guess will be that $P(x) := \displaystyle{-1\over f(x)f(-x)} + 2f_3$
is simply $\wp(x)$! Since $\wp$ satisfies $\wp'^2 = 4\wp^3
-g_2\wp -g_3$, by taking differentiation
one gets that $\wp''-6\wp^2$ is a constant $-g_2/2$. Thus we would like to
compute $P''-6P^2$. To simplify the presentation, let $g(x) = 1/f(x)$ and so
$P(x) = - g(x)g(-x) + 2f_3$. Then DE-1 and DE-2 takes the form DE-1* and DE-2*:
\begin{eqnarray*}
g(x)^2g(-x)^2 -a_1(g'(x)g(-x) \!\!\!&+&\!\!\! g'(-x)g(x))\\
 &+&\!\!\! g'(x)g'(-x) + g''(x)g(-x)=0 \\
6f_3g(x)g(-x) -2a_1(g'(x)g(-x) \!\!\!&+&\!\!\! g'(-x)g(x)) \\
 &+&\!\!\! 2g'(x)g'(-x) + g''(x)g(-x)=0.
\end{eqnarray*}
Since the only term that is not symmetric with respect to
$x\to -x$ is the last term $g''(x)g(-x)$, we must have
$$
g''(x)g(-x)=g''(-x)g(x).
$$
But then $(g'(x)g(-x) + g'(-x)g(x))'=0$ because it equals
$$
g''(x)g(-x) -g'(x)g'(-x) -g''(-x)g(x) + g'(-x)g'(x)=0.
$$
That is, $g'(x)g(-x) + g'(-x)g(x)$ is a constant. By expanding out the power series
one sees that this constant is $6f_4$. With these understood, then
\begin{eqnarray*}
P''(x)-6P(x)^2 \!\!\!&=&\!\!\! (-g'(x)g(-x)+g(x)g'(-x))'-6(- g(x)g(-x) + 2f_3)^2\\
  &=& \!\!\!-g''(x)g(-x) + 2g'(x)g'(-x) -g(x)g''(-x) \\
  && \!\!\!-6g^2(x)g^2(-x) + 24f_3g(x)g(-x) -24f_3^2.
\end{eqnarray*}
This is exactly $-6(\hbox{\rm DE-1*}) + 4(\hbox{\rm DE-2*}) + 12a_1f_4 -24f_3^2$. So
$$
P''(x)-6P(x)^2 =12a_1f_4-24f_3^2,
$$
which integrates into the Weierstrass equation with $g_2=-24a_1f_4+48f_3^2$. That is,
there exists periods lattice $\Lambda$ such that $P(x) = \wp(x)$. (When
the cubic curve is singular, $\Lambda$ is of  rank one and $\wp$ is a trigonometric function.)\par

In order to determine $f(x)$, recall that
$$
g'(x)g(-x) + g'(-x)g(x) = 6f_4.
$$
Also from the derivative of the equation $-g(x)g(-x)=\wp(x)-2f_3$ one gets
$$
-g'(x)g(-x) + g'(-x)g(x) = \wp'(x).
$$
This gives that $g'(x)g(-x)=(6f_4-\wp'(x))/2$. Hence
$$
{f'(x)\over f(x)}= -{g'(x)\over g(x)}=-{1\over 2}
{\wp'(x)-6f_4\over \wp(x)-2f_3}.
$$
Choose $z$ such that $\wp(z)=2f_3$
then $-g(z)g(-z)=\wp(z)-2f_3=0$. The choice of $z$ is up to sign, we choose the
one such that $g(-z)=0$. Then $(6f_4-\wp'(z))/2 =g'(z)g(-z)=0$, that is,
$6f_4=\wp'(z)$. So
$$
{f'(x)\over f(x)}=-{1\over 2}
{\wp'(x)-\wp'(z)\over \wp(x)-\wp(z)}.
$$
Let $a:=2f_3=\wp(z)$, $b:=3f_4=\wp'(z)$ and $g_2$ be
the corresponding coefficient in the Weierstrass equation.
Together with the fact that the extra factor $e^{kt}$ contributes simply an additive constant $k$ to
$(\log f(x))'$. We see that $S_2\subset S_3$.\hfill QED\par\medskip

\n{\bf Remark 3.1} (Coordinates of $S_3$)
(1) Since $S_3 = {\rm Spec}\,\mathbb{C}[k,a,b,g_2]$, $(k,a,b,g_2)$ is the algebraic
coordinates of $S_3$. In this paper we have ignored the integral structure of $S_3$ completely.
Readers interested in it may consult \cite{Totaro} for more details.
(2) When we represent $\Lambda=\mathbb{Z}\omega_1+\mathbb{Z}\omega_2$ and
$(a,b)=(\wp(z),\wp'(z))$, we obtain the analytic parameter system $(k,\omega_1,\omega_2,z)$.
This will be useful in the proofs of Theorem C and D.
(3) The obvious scaling $Q(sx)$ for $s\in \mathbb{C}^\times$ all correspond to
proportional genera (chern numbers). They correspond to recalling of the lattice,
so the complex elliptic genera is also usually regarded as depending on three parameters
$\tau=\omega_2/\omega_1$, $z$ and $k$ only.\\

\centerline{\bf \S4 Complex Elliptic Genera under Blowing-Up}\medskip

In this section we will complete the proof of Theorem B by showing that under the analytic
parameter system, $S_2$ contains precisely those points with $z$ not a 2 torsion point and
$S_1$ contains all points with $z$ a non-torsion point.\par

There are several equivalent definitions of complex elliptic genera in the literature.
It depends on the choices of elliptic-like functions and the parameter systems. It turns out that
our definition is very close to Krichever's \cite{Krichever}. Namely the complex
genus $\varphi$ defined by the four parameter ($k$, $\omega_1$, $\omega_2$, and $z$)
power series
$$
f(x) = e^{kx}\,e^{-\zeta(z)x}\,{\sigma_{\omega_1,\omega_2}(x)\sigma_{\omega_1,\omega_2}(-z)
\over \sigma_{\omega_1,\omega_2}(x - z)}.
$$

To see this, one may write out the function $f(x)$ we get in last section in terms of the
Weierstrass $\sigma$ function. Recall that $\zeta(x) = -\int^x \wp = 1/x + \cdots$ and
$\sigma(x) = e^{\int^x \zeta} = x +\cdots$, both are odd functions. Then we
have the well-known formula (see eg.\ \cite{Ahlfors})
$$
-{1\over 2}{\wp'(x)-\wp'(z)\over \wp(x)-\wp(z)}=\zeta(x) + \zeta(z) -\zeta(x+z).
$$
So $\log f(x) = \log\sigma(x) -\log\sigma(x+z) +\zeta(z)x + \lambda$. Where
$\lambda$ is easily seen to be $\log\sigma(z)$ by comparing coefficients. Since
the solution is always up to a normalization factor $e^{kx}$, the general solution
is thus given by
$$
f(x)=e^{kx}e^{\zeta(z)x}\,{\sigma(x)\sigma(z)\over\sigma(x+z)}.
$$
This agrees with Krichever's definition when we replace $z$ by $-z$. We will use our definition
throughout this paper. \par

Recall that $\sigma(z + \omega_i) = -e^{\eta_i(z + \omega_i/2)}\sigma(z)$ with
$\eta_i=\zeta(\omega_i/2)$. If $\lambda\in \mathbb{Z}$ then for $\vartheta(z) := \sigma(\lambda z + a)$,
this quasi-periodicity of $\sigma$ implies that
$$
\vartheta(z + \omega_i) = (-1)^\lambda e^{\eta_i(\lambda^2 z + \lambda^2 \omega_i/2 +\lambda a)}
\vartheta(z).
$$
Hence the following well-known fact\par\medskip

\n{\bf Lemma 4.1} {\it The function $\displaystyle{\prod\nolimits_{j = 1}^r {\sigma(\lambda_j z - a_j)
\over \sigma(\mu_j z - b_j)}}$ is elliptic, that is, doubly periodic if $\sum\lambda_j^2 = \sum\mu_j^2$,
$\sum\lambda_ja_j = \sum\mu_jb_j$ and $\sum \lambda_j \equiv \sum\mu_j \pmod{2}$.}\par\medskip

This Lemma also holds if $a_j$ and $b_j$ takes values in the nilpotent elements of some
commutative algebra. For example, even cohomologies or the Chow rings. In that case, the pole
of $1/\sigma(\lambda z -a)$ is still at $z = 0$ according to our definition (see the end of \S2).\par\medskip

Now we are ready to prove that $S_1$ (resp.\ $S_2$) contains all $f$'s in $S_3$ with $z$ a
non-torsion (resp.\ non $2$-torsion) point. In fact we will show that for $\varphi$ the complex elliptic
genera as defined above, the residue term in Theorem A for a blowing-up along codimension $r$
center is zero for $z$ not an $r$ torsion point. Notice that when the Weierstrass equations
define singular cubic curves, the periods lattice degenerates to rank one and
$\wp(x)$, $\zeta(x)$ and $\sigma(x)$ are all trigonometric
functions with the same defining properties as in the non-singular case. \par

Direct substitution shows that the residue is given by
$$
{\rm Res}_{\,t=0}\Big(e^{-kc_1(N)}e^{(r-1)(k+\zeta(z))t}{\sigma(t+z)\over \sigma(t)\sigma(z)}
\prod\nolimits_{i=1}^r{\sigma(n_i-t+z) \over \sigma(n_i-t)\sigma(z)}\,A(t,r)\Big).
$$
In order for this to be zero, by Lemma 4.1, if we choose (notice that $rz\not\in\Lambda$)
$$
A(t,r) = e^{-(r-1)(k+\zeta(z))t}\,{\sigma(t + rz)\sigma(z)\over\sigma(t+z)\sigma(rz)},
$$
then since $-rz + \sum_{i=1}^r(n_i + z) - z = \sum_{i=1}^r n_i - z$, we conclude that the above
power series is an elliptic function (with value in the Chow ring) and with $t=0$
the only pole (notice that the factor $\sigma(t+z)$ is canceled out). The contour
integration over a parallelogram domain now shows that the coefficient of degree $-1$ term
(the residue) must be identically zero. When the lattice degenerates to rank one, we
may use contour integral along the boundary of a thin tube and then take limits to conclude
the same result. Hence the proof.\par

It remains to consider the case that $z$ is a 2 torsion point. Suppose that $f(x)$ is a solution
to the functional equation. Let $z=\omega/2$ for some period $\omega$ and let
$e=\wp(\omega/2)$. As before we may also first assume that $k=0$. In this case
$f(x)=1/\sqrt{\wp(x)-e}$ is an odd function (this is the real elliptic genera considered
in \cite{HBJ}). Since $f(-x)=-f(x)$, the formula for $A(x)$ in \S3 reduces to $A(x)=a_1f(x) + f'(x)$. Plug
in this into the functional equation and replace $x$ by $-x$. After simplification we get
$$
f(x+y)=f'(x)f(y) + f'(y)f(x).
$$
By expanding out the power series and equating the coefficients term by term, one sees that
the solution, if exists, is uniquely determined by $f_3$. It is then easy to see that the general solution is
given by $f(x)=\sinh(sx)/s$ with $s^2=f'''(0)=6f_3\in\mathbb{C}$.
In this case the corresponding lattice is degenerate. The proof of Theorem B is completed.\hfill QED\par\medskip

\n{\bf Remark 4.2} (Uniqueness of $A(t,r)$) (1) At least when the cubic curve is smooth, we expect
that $A(t,r)$ is uniquely determined by $f$ and its existence is equivalent to that $z$ is not an $r$
torsion point. However, the author do not know a proof of this for $r\ge 4$. (2) Instead of
being an universal Jacobian factor, in specific cases
with $Y\to X$ fixed, if we allow $A(t,r)$ to have coefficients in cohomology classes, then the
choices is no longer unique. In fact any power series $e^{(r-1)kt}B(r,t)$ with the same
value as the chosen one at $t=0$ and with $B(r,t)$ satisfying the same
transformation property will do the job. We will see this non-uniqueness during the proof of
Theorem C.\par\medskip

\n{\bf Remark 4.3} There is an alternative way to prove that $S_2\subset S_3$ based
on the knowledge that $S_2$ contains at least those $f$'s of $S_3$ with $z$ not a 2 torsion.
The strategy is to compute the degree of freedom of $S_2$. Using the twisting $e^{kx}f(x)$,
we may first normalize $f(x)$ such that $f_2=0$. Now we expand out the functional equation via
power series in $x$ and $y$. Given $d\ge 2$ and $1\le p\le d-1$, comparing the coefficient of
$x^p y^{d-p}$ gives
\begin{eqnarray*}
\sum_{i=1}^{p-1}(-1)^i f_if_{d-i}C^d_p
  &=& \sum_{i=1}^p\sum_{j=1}^{d-j}C^{d-i-j}_{p-i}\big(a_if_j + (-1)^{d-i-j}a_jf_i\big)f_{d-i-j}\\
  && + \sum_{j=1}^{d-p} (-1)^{d-j}f_jf_{d-j}C^{d-j}_p + \sum_{i=1}^p(-1)^if_if_{d-i}C^{d-i}_{p-i}.
\end{eqnarray*}
For $d =2,3$, these are trivial identities. For $d =4$, all three equations
reduces to $a_2  - 3f_3 = 0$. So we may allow $f_3$ and $a_1$ to be arbitrary and then
solve $a_2 =  3f_3$. With this, for $d=5$, all the four equations are equivalent
and we get that $a_3 =  2f_4 + a_1f_3$ with $f_4$ arbitrary. For $d=6$, we get
$$
a_4=2a_1f_4 + {3 \over 2}f_3^2,\qquad f_5={3\over 10}f_3^2 + {3\over 5}a_1f_4.
$$

We want to solve $f(x)$ and $A(x)$ inductively. For $d\ge 6$, each time there appears two
new coefficients $a_{d-2}$ and $f_{d-1}$ and with $d-1$ relations. If we show that there are at
least two independent relations among those $d-1$'s then $f(x)$ and $A(x)$
are uniquely determined by $f_3$, $f_4$ and $a_1$, if they exist.\par

First of all, $a_{d-2}$ occurs only for $p=2 $ or $p=d-2$ (which by the symmetry of the functional
equation correspond to the same relation) as the term $1\cdot a_{d-2}f_1f_1 = a_{d-2}$.
But for all $p$, $f_{d-1}$ appears in the relations with coefficients
$$
C^d_p + (-1)^dC^d_p + (-1)^dC^{d-1}_p + C^{d-1}_{p-1},
$$
which is always nonzero. Hence there are at least two
independent relations and $f(x)$ is uniquely determined by $a_1$,
$f_3$ and $f_4$. Moreover if $f_4\ne 0$ then $f(x)$ is uniquely determined by
$f_3$, $f_4$ and $f_5$. By writing out $f(x)$ which defines the complex elliptic genera
with $k=0$, we find that
$$
f(x)=x + {a\over 2}x^3 + {b\over 6}x^4 + \Big({3a^2\over 8}-{g_2\over 40}\Big)x^5 +\cdots.
$$
This establishes an one to one correspondence between $S_2$ with $f_2=0$, $f_4\ne0$
and $S_3$ with $k=0$, $b\ne0$ --- that is with $z$ a non 2 torsion point.\par

In fact this research started from solving the functional equation inductively. We are confident
with our approach after it has been verified in Maple V and Mathematica up to degree 20
that $f(x)$ and $A(x)$ are uniquely solvable in $f_2$, $f_3$, $f_4$ and $a_1$.\\

\centerline{\bf \S5 The Change of Variable Formula}\medskip

\n{\bf Theorem 5.1} (Transition Formula) {\it Let $\varphi$ be the complex elliptic genera.
Let $E_i'$, $i=1,\ldots,p$ be $p$ irreducible divisors in $X$ and $e_i\in\mathbb{R}\backslash\{-1\}$.
Consider a blowing-up $\phi:Y\to X$ along smooth center of codimension $r$ with
exceptional divisor $E_0$. Let $\phi^*E_i' = E_i + m_iE_0$ with $E_i$ the proper transform of $E_i'$.
Let $D$ be an cycle in $X$. If $e_0:=\sum\nolimits_{i=1}^pe_im_i + (r-1)\ne -1$ then
$$
\int_D \prod\nolimits_{i=1}^p A(E_i', e_i + 1)\,K_\varphi(c(T_{X}))=
\int_{\phi^*D} \prod\nolimits_{i=0}^p A(E_i,e_i + 1)\,K_\varphi(c(T_Y)).
$$}

\n{\it Proof.} For the complex elliptic genera, Theorem A implies that for any power series
$F({\bf t})=F(t_1,\ldots,t_p)$ and cycles $D$, ${\bf E'}=(E_1',\ldots,E_m')$ in $X$,
$$
\int_D F({\bf E}')\,K_Q(c(T_X))=\int_{\phi^*D} F(\phi^*{\bf E}')\,A(E_0,r)\,K_Q(c(T_Y)).
$$
(Since $\bar\phi_*R_1(e) = 0$.) With this, the left hand side in the theorem becomes
$$
\int_Y \prod\nolimits_{i=1}^p A(\phi^*E_i',e_i + 1)\,A(E_0,r)\,K_\varphi(c(T_Y)).
$$
And the right hand side can be written as
$$
\int_Y \prod\nolimits_{i=1}^p A(\phi^*E_i'-m_iE_0,e_i + 1)\,A(E_0,e_0 + 1)\,K_\varphi(c(T_Y)).
$$

Now we plug in $A(t,r) = e^{-(r-1)(k+\zeta(z))t}
\displaystyle{\sigma(t +rz)\sigma(z)\over\sigma(t+z)\sigma(rz)}$ and analyze the
map $\phi$. The dominant variable $E_0$ is again replaced by the variable $t$ in the fiber
integration calculation.

The extra {\it Jacobian factors} of both integrals have the same exponential factor
$$
e^{-e_0(k+\zeta(z))t -\sum (-m_i)e_i(k+\zeta(z))t} = e^{-(e_0 -\sum m_ie_i)(k+\zeta(z))t}
= e^{-(r-1)(k+\zeta(z))t},
$$
and also the same
relevant transformation factor: for the first one, it is $e^{2\pi i(-rz +z)} = e^{-2\pi i(r-1)z}$;
for the second integral, the exponent is $2\pi i$ times
\begin{eqnarray*}
&&\!\!\!\sum\nolimits_{i=1}^p\big(m_i\phi^*E_i'+m_i(e_i + 1)z\big)-
  \sum\nolimits_{i=1}^p\big(m_i\phi^*E_i'+m_iz\big)-(e_0 +1 )z + z \\
&&\!\!\!= -(e_0 -\sum\nolimits_{i=1}^p m_i e_i)z = -(r-1)z.
\end{eqnarray*}
As in Remark 4.2, since both Jacobian factors become equal if we formally set $E_0=0$,
this implies that both have the same effect in the fiber integration computation of
$\bar\phi_*$ (which is zero), so both integrals are equal. \hfill QED\par\medskip

Now we prove Theorem C. Let $\phi': X'\to X$ be a birational morphism with
$K_{X'} = \phi'^*K_X + \sum_{i=1}^p e_iE_i'$. Consider a further blowing-up $\psi:Y\to X'$ along
a smooth center of codimension $r$ with $K_Y = \psi^*K_X + (r-1)E_0$. Let $\phi = \psi\circ\phi'$
and let $\psi^*E_i' = E_i + m_iE_0$. Then the canonical bundles satisfy the following relations
\begin{eqnarray*}
K_Y\!\!\!&=&\!\!\! \psi^*\Big(\phi'^*K_X + \sum\nolimits_{i=1}^p e_iE_i'\Big) + (r-1)E_0\\
  &=&\!\!\! \phi^*K_X + \sum\nolimits_{i=1}^p e_iE_i + \Big(\sum\nolimits_{i=1}^pe_im_i + (r-1)\Big)E_0.
\end{eqnarray*}
By applying Theorem 5.1 to the blowing-up $\psi:Y\to X'$ we conclude that
$$
\int_{X'} \prod\nolimits_{i=1}^p A(E_i', e_i + 1)\,K_\varphi(c(T_{X'}))=
\int_Y \prod\nolimits_{i=0}^p A(E_i,e_i + 1)\,K_\varphi(c(T_Y)).
$$
In particular, this proves Theorem C in the case that $\phi:Y\to X$ is a composite of blowing-ups
along smooth centers.\par

To prove Theorem C for general birational morphism $\phi$, we need to assume that $k=\mathbb{C}$
and make use of a recent result due to Wlodarsczyk and his co-workers \cite{AKMW} \cite{Wlodarsczyk},
namely the weak factorization theorem.
It says that (in characteristic zero) any birational map $f:X\cdots\!\!\to X'$ can be
factorized into composite of $f_i: X_i\cdots\!\!\to X_{i+1}$, $i=0,\ldots,q$ such that
$X_0 = X$, $X_{q+1}=X'$ and each $f_i$ is either a blowing-up or a blowing-down along smooth
center. We apply it to the morphism $\phi:Y\to X$.\par

Since the coefficient $e_i$ in front of $E_i$ is independent of the birational
model we choose, as long as the divisor $E_i$ has a nontrivial proper transform in
that model, they must transform correctly in all $f_i$. Theorem C then follows from the blowing-up case.
\hfill QED\\

\centerline{\bf \S6 $K$-equivalence Relation, Proof of Theorem D}
\centerline{\bf and the Main Conjectures}\medskip

Let $X$ be an $n$ dimensional complex normal $\mathbb{Q}$-Gorenstein variety.
Recall that $X$ has (at most) terminal (resp.\ canonical, resp.\ log-terminal) singularities
if there is a (hence for any) resolution $\phi:Y\to X$ such that in the canonical bundle relation
$K_Y=_\mathbb{Q} \phi^*K_X+\sum a_iE_i$,
we have that $a_i>0$ (resp.\ $a_i\ge 0$, resp.\ $a_i>-1$) for all $i$.
Here, the $E_i$'s vary among the prime components
of all the exceptional divisors. For two $\mathbb{Q}$-Gorenstein varieties $X$ and $X'$,
we say that $X$ and $X'$ are $K$-equivalent, written as
$X =_K X'$, if there is a smooth variety $Y$ and a birational correspondence
$(\phi,\phi'): X\leftarrow Y \to X'$, such that $\phi^*K_X =_\mathbb{Q} \phi'^*K_{X'}$.
Notice that this property does not depend on the choices of $Y$.\par

To get a feeling on the objects involved, let us recall
the following typical situations that lead to $K$-equivalence.
By definition, any composite of flops induces $K$-equivalence.
More generally, let $f\colon X\cdot\!\cdot\!\!\to X'$
be a birational map between two varieties with at most canonical
singularities such that $K_X$ (resp. $K_{X'}$) is nef along the
exceptional locus $Z\subset X$ (resp. $Z'\subset X'$), then $X=_K X'$.
In particular, this applies to birational minimal models \cite{Kollar} \cite{Wang}.
Also all cohomologically small resolutions of a singular variety, if they exist,
are all $K$-equivalent \cite{Totaro}.\par\medskip

We now prove Theorem D. It is clear that if
$X$ and $X'$ are $K$-equivalent proper smooth complex algebraic varieties, then
by the change of variable formula (Theorem C), we know that they have the same complex
elliptic genera at least for the parameter $z$ not an $r$ torsion point for $2\le r\le \dim\,X$.
But once we know that the complex elliptic genera coincide for generic $z$, they must coincide
by continuity (or specialization). Conversely, if a complex genus $\varphi$ is invariant
under $K$-equivalence then it is invariant under classical flops, hence by Totaro's theorem
\cite{Totaro} it must belong to the complex elliptic genera.\par

Now let $X$ be a complex $\mathbb{Q}$-Gorenstein variety with at most log-terminal singularities.
Take any resolution of singularities $\phi:Y\to X$ with $K_Y=_\mathbb{Q} \phi^*K_X+\sum e_iE_i$.
Since $e_i > -1$, one may simply define its complex elliptic genera to be
$$
\int_Y \prod A(E_i,e_i + 1)\,K_\varphi(c(T_Y)),
$$
where $A(t,r)$ is the same as before though now we plug in the variable $r$ by rational numbers.
Again this definition will cause difficulties for certain torsion values $z$, we avoid this
problem by using the universal complex elliptic genera instead of its various specializations.
In order to show that it is independent of the smooth model $Y$, suppose that $Y'\to X$ is another
resolution, then by using the weak factorization theorem, $Y$ and $Y'$ are connected
through blowing-ups and blowing-downs. Then Theorem 5.1 implies this independence
because the coefficient $e_i$ in front of $E_i$ is independent of the birational
model we choose.\par

Finally, it follows from Theorem B or Totaro's Theorem that there are no other genera which could be
defined on singular varieties such that they are compatible with the change of variable formula.\hfill QED\par\medskip

As in the case of birational minimal models,
we expect that any $K$-equivalence can be decomposed into composite of some {\it nice flops}.
However, this {\it rigid decomposition} is too hard to achieve at this moment. Instead,
we would like to state a series of conjectures on $K$-equivalent varieties, with the hope to
reduce the necessity of a rigid decomposition result for most potential applications.\par\medskip

\n{\bf Main conjectures on $K$-equivalence relation} ---
Fix a birational map $f:X\cdot\!\cdot\!\to X'$ between two
proper smooth complex varieties and let $T := \phi'_*\circ\phi^*$ be the cohomology correspondence
induced from a birational correspondence $(\phi,\phi'): X\leftarrow Y \to X'$ which extends $f$
and with smooth $Y$. $T$ is determined by the closure of the graph
$\bar\Gamma_f\subset X\times X'$ through the K\"unneth formula hence is independent of
the choice of $Y$. Suppose that $X=_K X'$.\par\smallskip

\n{\bf I} (canonical isomorphism) $T$ induces a canonical isomorphism on cohomologies,
which respects the rational Hodge structures:
$$
T:H^i(X,\mathbb{Q})\mathop{\longrightarrow}^\sim H^i(X',\mathbb{Q}).
$$

\n{\bf II} (quantum cohomology/K\"ahler moduli) Under part I,
$T$ also induces an isomorphism on the (big) quantum cohomology rings in
the sense of analytic continuations over the extended K\"ahler moduli spaces
(compare with \cite{Ruan}).\par\smallskip

\n{\bf III} (birational complex moduli) $X$ and $X'$ have canonically isomorphic (at least local)
complex moduli spaces. Moreover, suitably compactified polarized moduli spaces
should again be $K$-equivalent.\par\smallskip

\n{\bf IV} (soft decomposition) $X$ and $X'$ admit
symplectic deformations such that the $K$-equivalence relation deformed into copies of
classical flops.\par\medskip

Most of these conjectures are known in dimension three based on classification theoretic
results on flops in the minimal model theory \cite{Kollar} \cite{Ruan} \cite{KoMo}. Yet, the
techniques involved are unlikely to work in higher dimensions. It seems that IV will play a
key role toward the understanding of I, II and III. \par\medskip

\n{\bf Topological evidence for conjecture IV} --- Let $\Omega^U$ be the cobordism ring of
stably almost complex manifolds. For any $\mathbb{Q}$
algebra $R$, an $R$-valued complex genus defined in the topological way is a ring homomorphism
$\varphi:\Omega^U_\mathbb{Q}\to R$. A theorem due to Milnor says that the rational cobordism class
is determined exactly by all the chern numbers of the stable tangent bundle, or equivalently,
determined by all its complex genera. So in fact the topological definition of genera is the same as the
previous algebraic one. In terms of the cobordism theory, we may rephrase
Theorem D in the following way:\par

Totaro proved that `complex cobordism ring modulo classical flops' $=$
`complex elliptic genera'. Theorem D generalizes this to `complex cobordism ring
modulo $K$-equivalence' $=$ `complex elliptic genera'. That is, inside the complex cobordism
ring, the ideal generated by $X-X'$ for $X$ and $X'$ which are
related by classical flops are indeed the same as the seemingly much larger ideal generated by all
$X-X'$ where $X=_K X'$. So Conjecture IV is true up to complex cobordism.\\

\centerline{\bf \S7 Relations with Equivalence of Hodge Structures}\medskip

The proof of the equivalence of Hodge numbers sketched in \cite{Wang} uses the
theory of motivic integration developed by Denef and Loeser \cite{DeLo}. In fact, for $K$-equivalent
smooth complex projective varieties $X$ and $X'$ one has $[X]=[X']$ in a suitably completed
localized Grothendieck ring of algebraic varieties $\widehat{\cal M}$. As is remarked in \cite{DeLo},
the Hodge structure realization functor factors through this ring. Together with the fact that the
category of pure $\mathbb{Q}$-Hodge structures is semi-simple, we conclude that
$X$ and $X'$ have isomorphic $\mathbb{Q}$-Hodge structures on cohomologies.
However, this does not provide any canonical morphism between them.\par

Hodge numbers and Hodge structures determine a substantial part of the complex elliptic genera
and also give information to the complex moduli. For this, recall the formula in \cite{Totaro}:
$$
\varphi(X) = \chi\Big(X, K_X^{\otimes(-k)}\otimes\prod\nolimits_{m\ge 1}
(\Lambda_{-y^{-1}q^m}T\otimes\Lambda_{-y^{-1}q^{m-1}}T^*\otimes S_{q^m}T\otimes S_{q^m}T^*)\Big).
$$
Here we normalize the period lattice by $\omega_1=1$, $\omega_2=\tau$, also
$q = e^{2\pi i\tau}$, $y=e^{2\pi iz}$ and $T = T_X -n$ the rank zero virtual tangent bundle.
The twisted $\chi_y$-genus corresponds to the two parameter genera
$$
\chi_y(X) := \chi\Big(X,K_X^{\otimes(-k)}\otimes\Lambda_y T_X^*\Big),
$$
which is equivalent to knowing all $\chi(X,K_X^{\otimes(-k)}\otimes\Omega_X^p)$ for $p\ge 0$.
If $n=\dim X\le 11$, the twisted $\chi_y$ genus contains the same chern numbers as the
complex elliptic genera. So in this range, twisted $\chi_y$ genus contains precisely all
chern numbers that are invariant under the $K$-equivalence relation. If $n\le 4$, the twisted
$\chi_y$ genus contains all chern numbers, so all chern numbers are invariant under
$K$-equivalence for dimensions up to $4$.\par

It is clear that if $K_X$ is trivial, that is, $X$ is a Calabi-Yau manifold, then the twisted
$\chi_y$ genus becomes Hirzebruch's $\chi_y$ genus $\sum_{p\ge 0}\chi(X,\Omega_X^p)\,y^p$.
In particular, it is determined by the Hodge numbers. So the equivalence of elliptic genera
(that is, $k=0$) follows from the equivalence of Hodge numbers when $n\le 11$. But when
$n\ge 12$, the elliptic genera and Hodge numbers contain quite different type of
information.\par

For non Calabi-Yau manifolds, we can still use Hodge numbers to study twisted $\chi_y$
genus in some cases. First we show that:\par\medskip

\n{\bf Theorem 7.1} {\it Let $X$ and $X'$ be two $K$-equivalent smooth complex projective
varieties with $D\subset X$ and $D'\subset X'$ be base point free divisors
such that $\phi^*D = \phi'^*D'$ for some birational correspondence
$(\phi,\phi'): X\leftarrow Y \to X'$.
Then for all $\ell\in \mathbb{Z}$ and $p\ge 0$,
$\chi(X,\mathscr{O}(\ell D)\otimes\Omega^p)=\chi(X',\mathscr{O}(\ell D')\otimes\Omega^p)$.}\par\medskip

\n{\it Proof.} We use induction on dimension $n=\dim\,X=\dim\,X'$. This is trivial
if $n=1$, so we may assume that the theorem is true up to dimension $n-1\ge 1$.\par

By Bertini's theorem, we may assume that $D$ and $D'$ are smooth, irreducible and
corresponds to each other under proper transform.
Let $\tilde{D}$ be the proper transform of $D$ and $D'$ in $Y$ with
$\bar\phi := \phi|_{\tilde{D}}$ and $\bar\phi' := \phi'|_{\tilde{D}}$. Then
$\bar\phi^*K_D = (\phi^*(K_X + D))|_{\tilde{D}}=
(\phi'^*(K_{X'} + D'))|_{\tilde{D}}=\bar\phi'^*K_{D'}$. That is, $D$ and $D'$ are again
$K$-equivalent.\par

We will prove by induction on $\ell\in\mathbb{N}\cup\{0\}$ that
$\chi(X,\mathscr{O}(\ell D)\otimes\Omega^p)=\chi(X',\mathscr{O}(\ell D')\otimes\Omega^p)$,
which is enough since they are polynomials in $\ell$.
For $\ell = 0$ this is true by equivalence of Hodge numbers. So let $\ell \ge 1$. From
$$
0\to \mathscr{O}((\ell-1)D)\otimes\Omega^p\to\mathscr{O}(\ell D)\otimes\Omega^p
\to \Omega^p|_D\to 0,
$$
we get that
$$\chi(X,\mathscr{O}(\ell D)\otimes\Omega^p)=
\chi(X,\mathscr{O}((\ell-1)D)\otimes\Omega^p)+\chi(X,\mathscr{O}(\ell D)\otimes\Omega^p|_D).
$$
By the induction hypothesis on $\ell$, we only need to take care of the last term.
From $0\to T_D\to T_X|_D\to N_D \cong \mathscr{O}_D(D)\to 0$, we have that
$0\to \mathscr{O}_D(-D)\to \Omega^1|_D\to \Omega^1_D\to 0$, so $\Omega^p|_D =
\Omega^p_D\oplus\mathscr{O}_D(-D)\otimes \Omega^{p-1}_D$. Hence that
$$
\chi(X,\mathscr{O}(\ell D)\otimes\Omega^p|_D)=\chi(D,\mathscr{O}_D(\ell D)\otimes\Omega^p_D)+
     \chi(D,\mathscr{O}_D((\ell-1) D)\otimes\Omega^{p-1}_D).
$$
(For $p=0$, it is understood that the third term is 0.)
Since now $\dim\,D=\dim\,D'=n-1$, $D=_K D'$ and $\bar\phi^*(D|_D)=\bar\phi'^*(D'|_{D'})$,
the induction hypothesis on $n$ then concludes that
$\chi(X,\mathscr{O}(\ell D)\otimes\Omega^p|_D)=\chi(X',\mathscr{O}(\ell D')\otimes\Omega^p|_{D'})$.
This completes the proof. \hfill QED\par\medskip

\n{\bf Corollary 7.2} {\it Let $X$ and $X'$ be two smooth complex projective varieties which
are birational good minimal models, that is both $X$ and $X'$ have $K^{\otimes r}$ to be
base point free for some $r\in\mathbb{N}$. Then for all $\ell\in \mathbb{Z}$ and $p\ge 0$,
$$
\chi(X,K_X^{\otimes \ell}\otimes\Omega^p)=\chi(X',K_{X'}^{\otimes \ell}\otimes\Omega^p).
$$}

\n{\it Proof.} Simply take $D=K_X^{\otimes r}$ and $D'=K_{X'}^{\otimes r}$ in the above
theorem and notice that the equality holds for all $\ell\in r\mathbb{N}$ implies that
it holds for all $\ell\in\mathbb{Z}$, since both terms are polynomials in $\ell$.
\hfill QED\par\medskip

Since birational minimal models are $K$-equivalent, Corollary 7.2 is just a special case of
Theorem C. But this alternative discussion has another aspect. Instead of using the
Euler characteristic functor, if we write out the two corresponding long exact sequences
for $X$ and $X'$ in the above proof, we may conclude inductively that
under conjecture I,
$$
H^q(X,K_X^{\otimes \ell}\otimes\Omega^p)\cong H^q(X',K_{X'}^{\otimes \ell}\otimes\Omega^p)
$$
for all $\ell\in r\mathbb{N}\cup \{0\}$.
It is likely that this will also hold for all $\ell\in\mathbb{Z}$. In that case we may
take $\ell=-1$ and use Serre duality theorem to get for all $i\ge 0$ that
$$
H^i(X,T_X)\cong H^i(X',T_{X'}).
$$
We hope that this will be useful in attacking Conjecture III concerning the
birational moduli spaces.\\

\end{document}